\documentclass[11pt]{article}
\usepackage{amssymb, amsmath, amsthm}

\usepackage{amsfonts,mathtools,mathrsfs,color}
\usepackage{bbm, dsfont}
\usepackage{algorithm}
\usepackage{algpseudocode}
\usepackage{graphicx}
\usepackage{xcolor}
\usepackage{tocloft}
\usepackage{bigints}
\usepackage{makecell}
\usepackage{centernot}

\usepackage{bm}
\usepackage{stmaryrd}

\usepackage[colorlinks]{hyperref}
\definecolor{coolblack}{rgb}{0.0, 0.0, 0.230}
\hypersetup{
  colorlinks,
  citecolor=coolblack,
  linkcolor=coolblack,
  urlcolor=coolblack}

\allowdisplaybreaks

\newtheorem{theorem}{Theorem}
\newtheorem*{theorem*}{Theorem}

\title{On the 3-adic Valuation of a Cubic Binomial Sum}
\author{Valentio Iverson\footnote{Department of Pure Mathematics, University of Waterloo, Waterloo, ON, Canada. \\ Email: \texttt{viverson@uwaterloo.ca}}}
\date{\today}

\begin{document}

\maketitle

\begin{abstract}
In this short note, we prove a conjecture recently posed by Alekseyev, Amdeberhan, Shallit, and Vukusic on the $3$-adic valuation of a cubic binomial sum.
\end{abstract}

\section{Introduction}
The study of $p$-adic valuations of combinatorial sequences and binomial sums is a classical and active area of research in number theory, often revealing that these arithmetic structures possess well-behaved $p$-adic valuations (see, for example, \cite{alekseyev2025padicvaluationsvalueslegendre, mikic2023divisibility}).

Recently, Alekseyev, Amdeberhan, Shallit, and Vukusic \cite{alekseyev2025padicvaluationsvalueslegendre} investigated the $p$-adic valuations of Legendre polynomials and related integer sequences. In their paper, they posed a specific conjecture regarding the $3$-adic valuation of the cubic binomial sum $\sum_{r = 0}^n \binom{n}{r}^3 2^r$. Specifically, they proposed an explicit formula for the $3$-adic valuation of this sum that depends strictly on the parity of $n$ and the base-$3$ sum of digits function, $s_3(\cdot)$.

In this short note, we resolve this conjecture in the affirmative. Our proof is elementary and direct. We proceed by applying MacMahon's identity to transform the cubic binomial sum into a form more amenable to $p$-adic analysis. Subsequently applying Legendre's formula to the transformed terms, we establish lower bounds on the $3$-adic valuations of the individual summands, employing techniques similar to those used in \cite{alekseyev2025padicvaluationsvalueslegendre}. This approach allows us to efficiently isolate the single dominating term that uniquely dictates the valuation of the entire sum.

Our main theorem is as follows:
\begin{theorem} 
For every integer $n \ge 1$, we have
\[
\nu_3 \left( \sum_{r = 0}^n \binom{n}{r}^3 2^r \right)
=
\begin{cases}
s_3 \left( \dfrac{n - 1}{2} \right) + 1 &\text{if } n \ \text{is odd}, \\[1ex]
s_3 \left( \dfrac{n}{2} \right) &\text{if } n \ \text{is even}.
\end{cases}
\]
\end{theorem} 

\section{Proof of Theorem 1}
Let 
\[ S_n := \sum_{r = 0}^n \binom{n}{r}^3 2^r.\]
We first note that for $n = 1$, we have $S_1 = 3$ from which it is easily verified that $\nu_3(S_1) = 1 =s_3(0) + 1$, as the theorem has suggested. Therefore, from here onwards, we shall prove the theorem for $n \ge 2$. Our starting point is the following classical identity:

\begin{theorem}[MacMahon's Identity \cite{macmahon1902sums}]
    For every positive integer $n$, we have
    \[ \sum_{k = 0}^n \binom{n}{k}^3 x^k y^{n - k} = \sum_{k = 0}^{\lfloor n/2 \rfloor} \binom{n}{2k} \binom{2k}{k} \binom{n + k}{k} x^k y^{k} (x + y)^{n - 2k}. \]
\end{theorem}
Substituting $x = 2$ and $y = 1$ into MacMahon's identity, we can rewrite $S_n$ as follows:
\begin{align*}
    S_n &= \sum_{r = 0}^n \binom{n}{r}^3 2^r \\ 
    &= \sum_{r=0}^{\lfloor n/2 \rfloor} \binom{n}{2r} \binom{2r}{r} \binom{n + r}{r} 2^r 3^{n - 2r}  \\ 
    &= \sum_{r = 0}^{\lfloor n/2 \rfloor} \binom{n + r}{3r} \binom{2r}{r} \binom{3r}{r} 2^r 3^{n - 2r}, 
\end{align*}
where the last equality can be checked by a straightforward computation. 

Let us denote the summands by 
\[ A_r = \binom{n + r}{3r} \binom{2r}{r} \binom{3r}{r} 2^r 3^{n - 2r}, \quad 0 \le r \le \lfloor n/2 \rfloor. \]
Recall that by Legendre's formula, for every positive integer $k$, we have 
\begin{equation} \label{eq:valuation}
\begin{split}
\nu_3 \left( \binom{2k}{k} \binom{3k}{k} \right) &= \nu_3 \left(\frac{(3k)!}{(k!)^3}\right) \\
&= \nu_3((3k)!) - 3 \nu_3(k!) \\
&= \frac{3k - s_3(3k)}{2} - 3 \cdot \frac{k - s_3(k)}{2} \\
&= s_3(k).
\end{split}
\end{equation}
Using this fact, we evaluate $S_n$ by considering the parity of $n$:
\begin{enumerate}
    \item If $n$ is even, we can write $n = 2m$ for some positive integer $m$. We aim to show that $\nu_3(S_{2m}) = s_3(m)$. Notice that 
    \begin{align*}
        S_{2m} &= \sum_{r = 0}^m \binom{2m + r}{3r} \binom{2r}{r} \binom{3r}{r} 2^r 3^{2m - 2r}. 
    \end{align*}
    The term for $r = m$ is given by $A_m = \binom{3m}{3m} \binom{2m}{m} \binom{3m}{m} 2^m$. Thus, by Legendre's formula and Equation \ref{eq:valuation}, we have
    \begin{align*} \nu_3(A_{m}) &= \nu_3 \left( \binom{2m}{m} \binom{3m}{m} 2^m  \right) =  s_3(m).
    \end{align*} 
    For every $0 \le r \le m - 1$, we observe the following lower bound: 
    \begin{align*}
        \nu_3(A_r) &= \nu_3 \left( \binom{2m + r}{3r} \binom{2r}{r} \binom{3r}{r} 2^r 3^{2m - 2r} \right) \\ 
        &\ge \nu_3 \left( \binom{2r}{r} \binom{3r}{r} \right) + 2(m - r) \\ 
        &= s_3(r) + 2(m- r) \\ 
        &\ge 1 + (m - r) + s_3(r) \\ 
        &\ge 1 + s_3(m - r) + s_3(r) \ge 1 + s_3(m),
    \end{align*}
    where the last inequality holds by the subadditivity of the sum of digits function.
    Because the valuation of $A_m$ is strictly less than the valuation of all other terms $A_r$, it dictates the valuation of the sum. Therefore,
    \[ \nu_3(S_{2m}) = \nu_3 \left( \sum_{r = 0}^m A_r \right) = s_3(m). \]
    
    \item If $n$ is odd, we can write $n = 2m + 1$ for some positive integer $m$. We have
    \[ S_{2m + 1} = \sum_{r = 0}^m \binom{2m + 1 + r}{3r} \binom{2r}{r} \binom{3r}{r} 2^r 3^{2m + 1 - 2r}. \]
    The term for $r = m$ is $A_m = \binom{3m + 1}{3m} \binom{2m}{m} \binom{3m}{m} 2^m 3$, yielding
    \begin{align*}
        \nu_3(A_m) &= \nu_3 \left( 2^m \cdot 3 \cdot (3m + 1) \binom{2m}{m} \binom{3m}{m} \right) \\ 
        &= 1 + \nu_3 \left( \binom{2m}{m} \binom{3m}{m} \right) = 1 + s_3(m). 
    \end{align*}
    Furthermore, for every $0 \le r \le m - 1$, we see that 
    \begin{align*}
        \nu_3(A_r) &= \nu_3 \left(\binom{2m + 1 + r}{3r} \binom{2r}{r} \binom{3r}{r} 2^r 3^{2m + 1 - 2r}  \right)  \\ 
        &\ge \nu_3 \left( \binom{2r}{r} \binom{3r}{r} \right) + (2m + 1 - 2r) \\ 
        &= s_3(r) + 2(m - r) + 1 \\ 
        &\ge 2 + (m - r) + s_3(r) \\ 
        &\ge 2 + s_3(m - r) + s_3(r) \ge 2 + s_3(m),
    \end{align*}
    similarly as before. Once again, the valuation is dictated by $A_m$, showing that
    \[ \nu_3(S_{2m + 1}) = \nu_3 \left( \sum_{r = 0}^m A_r \right) = s_3(m) + 1. \]
\end{enumerate}

This concludes the proof, showing that for every integer $n \ge 1$, 
\[ \nu_3 \left( \sum_{r = 0}^n \binom{n}{r}^3 2^r \right) = \begin{cases} s_3 \left( \frac{n - 1}{2} \right) + 1 &\text{if } n \ \text{is odd}, \\ s_3 \left( \frac{n}{2} \right) &\text{if } n \ \text{is even}. \end{cases} \]

\section{Acknowledgements}
The author would like to thank Ingrid Vukusic for proofreading and providing suggestions for improving the readability of the draft.

\end{document}